\documentclass[12pt,letterpaper]{article}
\usepackage{amsmath,amssymb,amsthm}
\usepackage[margin=1in]{geometry}
\usepackage{microtype}
\usepackage{url}
\usepackage{tikz}
\usepackage{enumitem}

\newtheorem{thm}{Theorem}[section]

\newtheorem{cor}[thm]{Corollary}
\newtheorem{lemma}[thm]{Lemma}

\newtheorem{claim}[thm]{Claim}

\theoremstyle{definition}

\theoremstyle{definition}

\theoremstyle{definition}
\newtheorem{defn}[thm]{Definition}

\theoremstyle{definition}

\theoremstyle{definition}

\theoremstyle{definition}

\theoremstyle{remark}

\usepackage[colorlinks=true, linkcolor=blue]{hyperref}

\newcommand{\beq}[1]{\begin{equation}\label{#1}}
\newcommand{\enq}[0]{\end{equation}}

\newcommand{\cE}[0]{{\cal E}}

\newcommand{\cG}[0]{{\cal G}}

\newcommand{\cI}[0]{{\cal I}}

\newcommand{\cO}[0]{{\cal O}}

\newcommand{\gD}[0]{\Delta}

\newcommand{\gS}[0]{\Sigma}

\begin{document}
\renewcommand{\thefootnote}{\fnsymbol{footnote}}

\title{Independent sets in the discrete hypercube}

\author{David Galvin\thanks{Department of Mathematics, University of Notre Dame, Notre Dame IN 46556; {\tt dgalvin1@nd.edu}. Research supported in part by the Simons foundation.}}

\maketitle

\begin{abstract}
In this expository note we describe a proof due to A. Sapozhenko that the number of independent
sets in the discrete $d$-dimensional hypercube $Q_d$ is
asymptotically $2 \sqrt{e} 2^{2^{d-1}}$ as $d$ tends to infinity.
\end{abstract}

\section{Introduction}

The focus of this note is the discrete hypercube $Q_d$. This is the
graph on vertex set $\{0,1\}^d$ with two strings adjacent if they
differ on exactly one coordinate. It is a $d$-regular bipartite
graph with bipartition classes $\cE$ and $\cO$, where $\cE$ is the
set of vertices with an even number of $1$'s. Note that
$|\cE|=|\cO|=2^{d-1}$. (For graph theory basics, see e.g.
\cite{MGT},~\cite{Diestel}).

An {\em independent set} in $Q_d$ is a set of vertices no two of
which are adjacent. Write $\cI(Q_d)$ for the set of independent sets
in $Q_d$. A trivial lower bound on $|\cI(Q_d)|$ is $2^{2^{d-1}}$:
each of the $2^{2^{d-1}}$ subsets of $\cE$ is an independent set. A
beautiful result of Korshunov and Sapozhenko
\cite{KorshunovSapozhenko} asserts that this trivial bound is not
far off the truth.

\begin{thm} \label{thm-sapkor}
$$
|\cI(Q_d)| = 2 \sqrt{e}(1+o(1))2^{2^{d-1}} ~~~\mbox{as} ~~~d
\rightarrow \infty.
$$
\end{thm}

The purpose of this expository note is to describe Sapozhenko's simplification
\cite{Sapozhenko2} of the proof of Theorem \ref{thm-sapkor}. The simplified proof of Theorem
\ref{thm-sapkor} depends on a technical lemma bounding the number of
subsets of $\cE$ of a given size whose neighbourhood in $\cO$ is of
a given size. This lemma, which appears in \cite{Sapozhenko}, is
stated in (close to) full generality in Section \ref{sec-saplemma}
and its proof is given in Section \ref{sec-saplemmaproof}. Section
\ref{sec-tools} establishes notation and gathers together in a
single place all the tools that we need. Mostly these are in the
form of references, but one important tool, a familiar isoperimetric
inequality in the hypercube, is proven in full, since we are not
aware of an explicit presentation of it in the literature in
English. The proof of Theorem \ref{thm-sapkor} appears in Section
\ref{sec-sapkorproof}.

\bigskip

To improve the trivial lower bound $|\cI(Q_d)| \geq 2^{2^{d-1}}$ to
that given by Theorem \ref{thm-sapkor}, we consider not just
independent sets which are confined purely to either $\cE$ or $\cO$.
It is easy to see that there are
$$
2^{d-1} 2^{2^{d-1}-d}=\frac{1}{2}2^{2^{d-1}}
$$
independent sets that have just one vertex from $\cO$, and more
generally approximately
$$
\frac{\left(\frac{1}{2}\right)^{k}}{k!} 2^{2^{d-1}}
$$
independent sets which have exactly $k$ non-nearby vertices from
$\cO$, for small $k$ (by ``non-nearby'' it is meant that there are
no common neighbours between pairs of the vertices). Indeed, there
are approximately ${2^{d-1} \choose k}$ ways to choose the $k$
vertices from $\cO$; these vertices together have a neighbourhood of
size $kd$, so there are $2^{2^{d-1}-kd}$ extensions of the $k$
vertices to an independent set. Summing over $k$ from $0$ to any
$\omega(d)$ we obtain the $\sqrt{e}$ in the lower bound.

To motivate the upper bound, consider what happens when we count
independent sets that have exactly two nearby vertices from $\cO$,
(i.e., two vertices with a common neighbour). There are
approximately $d^22^d$ choices for this pair (as opposed to
approximately $2^{2d}$ choices for a pair of vertices without a
common neighbour), since once the first vertex has been chosen the
second must come from the approximately $d^2$ vertices at distance
two from the first. There are approximately $2^{2^{d-1}-2d}$
extensions of the pair to an independent set (roughly the same as
the number of extensions in the case of the pair of vertices without
a common neighbour). The critical point here is that a pair of
vertices from $\cO$ has a neighbourhood size of approximately $2d$,
whether or not the vertices are nearby. This is because a pair of
vertices in $Q_d$ has at most two common neighbours. Thus we get an
additional contribution of approximately $2^{2^{d-1}-d}$ to the
count of independent sets from those sets with two nearby vertices
from $\cO$ (negligible compared to the addition contribution of
approximately $\frac{1}{4}2^{2^{d-1}}$ to the count from those sets
with two non-nearby vertices from $\cO$). The main work in the upper
bound is the correct extension of this observation to the
observation that the only non-negligible contribution to the count
is from independent sets that on one side consist of a set of
vertices with non-overlapping neighbourhoods. This in turn amounts
to showing that there is a negligible contribution from those
independent sets which are ``$2$-connected'' on one side (i.e., are
such that between any two vertices on one side, there is a path in
the cube every second vertex of which passes through the independent
set.) This, finally, entails proving the technical lemma bounding
the number of subsets of $\cE$ of a given size whose neighbourhood
in $\cO$ is of a given size.

\section{Notation and tools} \label{sec-tools}

\subsection{Notation}

Let $\Sigma$ be a $d$-regular bipartite graph on vertex set $V$ with
bipartition classes $X$ and $Y$. For $A \subseteq V$ we write $N(A)$
for the set of vertices outside $V$ that are neighbours of a vertex
in $A$, and $N(u)$ for $N(\{u\})$. We write $\rho(u,v)$ for the
length of the shortest $u$-$v$ path in $\Sigma$.

We define the {\em closure} of $A$ to be
$$
[A]=\{v \in V:N(v)\subseteq N(A)\}
$$
and say that $A$ is {\em closed} if $[A]=A$. We say that $A
\subseteq X$ is {\em small} if $[A] \leq |X|/2$.

We say that $A$ is {\em $k$-linked} if for every $u,v\in A$ there is
a sequence $u=u_0, u_1, \ldots, u_l=v$ in $A$ with
$\rho(u_i,u_{i+1})\leq k$ for $i = 0, \ldots, l-1$. Note that if $A$
is $2$-linked, then so is $[A]$. For any $k$, a set $A$ can be
decomposed into its maximal $k$-linked subsets; we refer to these as
the {\em $k$-components} of $A$.

For $u \in V$ and $A,B \subseteq V$ we write $\nabla(A)$ for the set
of edges having one end in $A$ and $\nabla(A,B)$ for the set of
edges having one end in each of $A, B$; $N_B(u)=N(u) \cap B$,
$N_B(A)=N(A) \cap B$ and $d_B(u)=|N_B(u)|$. Set $\rho(u,A)=\min_{w
\in A}\{\rho(u,w)\}$.

Given $A \subseteq X$ we always write $G$ for $N(A)$, $a$ for
$|[A]|$, $g$ for $|G|$ and set $t=g-a$.

Throughout we use $\log$ for the base $2$ logarithm. We do not track
constants in the proofs; all implied constants in $O$ and $\Omega$
notation are independent of $d$. We will always assume that $d$ is
sufficiently large to support our assertions.

\subsection{Tools}

The following easy lemma is from \cite{Sapozhenko}.

\begin{lemma} \label{nearbysets}
If $A$ is $k$-linked, and $T \subseteq V$ is such that
$\rho(u,T)\leq l$ for each $u \in A$ and $\rho(v,A)\leq l$ for each
$v \in T$, then $T$ is $(k+2l)$-linked.
\end{lemma}

We need a lemma that bounds the number of $k$-linked subsets of $V$.
The infinite $\Delta$-branching rooted tree contains precisely
$$
\frac{{\Delta n \choose n}}{(\Delta-1)n+1}
$$
rooted subtrees with $n$ vertices (see e.g. Exercise 11 (p. 396) of
\cite{Knuth}). This implies that the number of $n$-vertex subsets of
$V$ which contain a fixed vertex and induce a connected subgraph is
at most $(e\Delta)^{n}$. We will use the following easy corollary
which follows from the fact that a $k$-linked subset of $\Sigma$ is
connected in a graph with all degrees $O(d^{k+1})$.

\begin{lemma} \label{Tree}
For each fixed $k$, the number of $k$-linked subsets of $V$ of size
$n$ which contain a fixed vertex is at most $2^{O(n\log d)}$.
\end{lemma}

\medskip

The next lemma is a special case of a fundamental result due to
Lov\'asz \cite{Lovasz} and Stein \cite{Stein}. For a bipartite graph
$\Gamma$ with bipartition $P \cup Q$, we say that $Q' \subseteq Q$
{\em covers} $P$ if each $p \in P$ has a neighbour in $Q$.

\begin{lemma} \label{lovaszstein}
If $\Gamma$ as above satisfies $|N(p)| \geq a$ for each $p \in P$
and $|N(q)| \leq b$ for each $q \in Q$, then $P$ is covered by some
$Q' \subseteq Q$ with
$$|Q'| \leq (|Q|/a)(1 + \ln b).$$
\end{lemma}

We also use a result concerning the sums of binomial coefficients
which follows from the Chernoff bounds \cite{Ch} (see also
\cite{RG}, p.11):
\begin{equation} \label{inq-binomial}
\sum_{i=0}^{[cN]}{N\choose i} \leq 2^{H(c)N}~~~~~\mbox{for $c \leq
\frac{1}{2}$},
\end{equation}
where $H(x)=-x\log x -(1-x)\log(1-x)$ is the usual binary entropy
function and $[x]$ denotes the integer part of $x$. Also, throughout
we will use (usually without comment) a simple observation about
sums of binomial coefficients: if $k=o(n)$, we have
\begin{eqnarray}
\sum_{i \leq k} {n \choose i} & \leq & (1+O(k/n)){n \choose k} \nonumber \\
& \leq & (1+O(k/n))(en/k)^k \nonumber \\
& \leq & \exp_2\left\{(1+o(1))k\log(n/k)\right\}. \nonumber
\end{eqnarray}

\subsection{Isoperimetry in the cube}

A {\em Hamming ball centered at $x_0$} in $Q_d$ is any set of
vertices $B$ satisfying
$$
\{u \in V(Q_d)~:~\rho(u,x_0) \leq k\} \subseteq B \subseteq \{u \in
V(Q_d)~:~\rho(u,x_0) \leq k+1\}
$$
for some $k<d$, where $\rho$ is the usual graph distance (which in
this case coincides with the Hamming distance on $\{0,1\}^d$). An
{\em even} (resp. {\em odd}) {\em Hamming ball} is a set of vertices
of the form $B \cap {\cE}$ (resp. $B \cap {\cO}$) for some Hamming
ball $B$. We use the following result of K\"orner and Wei
\cite{KornerWei}. (A similar isoperimetric bound of Bezrukov
\cite{Bezrukov} would also suffice).

\begin{lemma} \label{kornerandwei}
For every $C \subseteq {\cE}$ (resp. ${\cO}$) and $D \subseteq
V(Q_d)$, there exists an even (resp. odd) Hamming ball $C'$ and a
set $D'$ such that $|C'|=|C|$, $|D'|=|D|$ and $\rho(C',D') \geq
\rho(C,D)$.
\end{lemma}

The following is a well-known isoperimetric inequality in $Q_d$
(see, e.g. \cite[Lemma 1.3]{KorshunovSapozhenko}).
\begin{claim} \label{clm-hyper.iso}
For $A \subseteq \cE$ (or $\cO$) with $|A| \leq M/2$ we have
$$
\frac{|N(A)|-|A|}{|N(A)|} = \Omega(1/\sqrt{d}).
$$
\end{claim}

\medskip

\noindent {\em Proof: }Without loss of generality, we may assume
that $A \subseteq {\cal E}$. Let such an $A$ be given. Applying
Lemma \ref{kornerandwei} with $C=A$ and $D=V(Q_d) \setminus (A \cup
N(A))$, we find that there exists an even Hamming ball $A^\prime$
with $|A^\prime|=|A|$ and $|N(A)| \geq |N(A^\prime)|$. So we may
assume that $A$ is an even Hamming ball.

We consider only the case where $A$ is centered at an even vertex,
without loss of generality $\underline{0} := \{0,\dots,0\}$, the
other case being similar. In this case,
$$
\{v \in {\cE}~:~\rho(v,\underline{0}) \leq k\} \subseteq A \subseteq
\{v \in {\cE}~:~\rho(v,\underline{0}) \leq k+2\}
$$
for some even $k \leq d/2$ (the bound on $k$ coming from the fact
that $|A|\leq M/2$). For each $0 \leq i \leq (k+2)/2$, set $B_i=A
\cap \{v: \rho(v,\underline{0})=2i\}$, and $N^+(B_i)=N(B_i)\cap\{u:
\rho(u,\underline{0})=2i+1\}$. It is clear that $N(A)=\cup_{0 \leq i
\leq (k+2)/2} N^+(B_i)$ and that
$$
\cup_{0 \leq i \leq k/2} N^+(B_i) = \{v \in
\cO~:~\rho(v,\underline{0}) \leq k+1\}.
$$
Also, observe that for all $i$
\begin{equation} \label{LYM}
\frac{|B_i|}{|N^+(B_i)|} \leq \frac{2i+1}{d-2i},
\end{equation}
from which it follows that
$$
|N^+(B_{(k+2)/2})|-|B_{(k+2)/2}| \geq
\frac{-10}{d-4}|N^+(B_{(k+2)/2})| \geq \frac{-20}{d}{d \choose k+3}.
$$
Indeed, (\ref{LYM}) is an equality except when $i=(k+2)/2$, in which
case it follows from the fact that each vertex in $B_{(k+2)/2}$ has
exactly $d-(k+2)$ neighbours in $N^+(B_{(k+2)/2})$, and each vertex
in $N^+(B_{(k+2)/2})$ has at most $(k+2)+1$ neighbours in
$B_{(k+2)/2}$.

We deal first with the case $k \leq d/4$. In this case, (\ref{LYM})
gives
$$
\frac{|B_i|}{|N^+(B_i)|} \leq \frac{2}{3}
$$
and so $(|N(A)|-|A|)/|N(A)| \geq 1/3$.

For $d/4 < k \leq d/2$ and $c$ any constant, we claim that
\begin{equation}\label{inq-binom_sum}
\sum_{i=0}^{k+c} {d \choose i} = O\left(\sqrt{d}{d \choose
k+c}\right),
\end{equation}
\noindent where the constant in the $O$ depends on $c$. For $k+c\leq
d/2$, this follows from (\ref{inq-binomial}):
$$
\sum_{i=0}^{k+c} {d \choose i} \leq 2^{H\left(\frac{k+c}{d}\right)d}
= O\left(\sqrt{d}{d \choose k+c}\right),
$$
the equality being an easy consequence of Stirling's approximation,
while for $k+c \geq d/2$ we have
$$
\sum_{i=0}^{k+c} {d \choose i} \leq 2^d = O\left(\sqrt{d}{d \choose
k+c}\right),
$$
once again from Stirling's approximation. For $d/4 < k \leq d/2$, we
therefore have
\begin{eqnarray}
{|N(A)|-|A|} & = & {|N^+(B_{(k+2)/2})|-|B_{(k+2)/2}| +
\sum_{i=0}^{k/2}
{d \choose 2i+1} - {d \choose 2i}} \nonumber \\
& \geq & {\frac{-20}{d}{d \choose k+3} + \sum_{i=0}^{k/2} {d-1
\choose 2i}\frac{d(d-1)}{(2i+1)(d-2i)}}
\nonumber \\
& \geq & \frac{8(d-1)}{3(d+2)}{d-1 \choose k}-\frac{20}{d}{d \choose
k+3}\nonumber \\
& = & \Omega\left({d \choose k+3}\right),\label{numer}
\end{eqnarray}
and
\begin{eqnarray}
{|N(A)|} & = & {|N^+(B_{(k+2)/2})| + \sum_{i=0}^{k/2} {d \choose
2i+1}}
\nonumber\\
& \le & {\sum_{i=0}^{k+3} {d \choose i}} =O\left({\sqrt{d}{d \choose
k+3}}\right),\label{denom} \
\end{eqnarray}
the last equality following from (\ref{inq-binom_sum}). Combining
(\ref{numer}) and (\ref{denom}), we may now conclude that
\[
\frac{|N(A)|-|A|}{|N(A)|} = \Omega\left(\frac{1}{\sqrt{d}}\right).
\]
\qed

Considerably stronger results can be obtained if we impose more
conditions on $|A|$. In that direction we need only the following
simple lemma.

\begin{lemma} \label{boundsondeltaAsmall}
For $A \subseteq \cE$ (or $\cO$),
$$
\mbox{if $|A|<d^{O(1)}$, then $|A| \leq O(1/d)|N(A)|$}.
$$
\end{lemma}

\medskip

\noindent {\em Proof: }If $|A|<d^{O(1)}$, then we have $k =O(1)$ in
the notation of Claim \ref{clm-hyper.iso}, and repeating the
argument of that lemma we get $|A| \leq O(1/d)|N(A)|$.
\qed

\section{Sapozhenko's graph lemma} \label{sec-saplemma}

Let $\gS$ be a $d$-regular, bipartite graph with bipartition classes
$X$ and $Y$. The {\em co-degree} of $\gS$ is the maximum of $|N(x)
\cap N(y)|$ over all pairs $(x,y) \in Y \times Y$. For each $a$ and
$g$ and $v \in \cO$, set
$$
{\cal G}(a,g,v) = \{A \subseteq X \mbox{ $2$-linked}:|[A]|=a,
~|G|=g~\mbox{and}~v \in G\}.
$$

\begin{lemma} \label{lem-mainsaplemma}
For each pair of constants $c,\gD_2>0$ there is a constant
$c'=c'(c,\gD_2)>0$ such that the following holds. If $\gS$ is a
$d$-regular bipartite graph with partition classes $X$ and $Y$ and
with co-degree $\gD_2$, then for $g \geq d^4$ and
\begin{equation} \label{lowerboundont}
t > \frac{c g\log^3 d}{d^2}
\end{equation}
we have
$$
|{\cal G}(a,g,v)| \leq 2^{g-\frac{c't}{\log d}}.
$$
\end{lemma}

\section{Proof of Theorem \ref{thm-sapkor}}
\label{sec-sapkorproof}

The key observation is the following consequence of Lemma
\ref{lem-mainsaplemma}.

\begin{cor} \label{cor-mainsaplemma}
$$
\sum_{\mbox{$A \subseteq \cE$, $A$ small}} 2^{-|N(A)|} =
(1+o(1))\sqrt{e} ~~~~\mbox{as $d \rightarrow \infty$}.
$$
\end{cor}

\medskip \noindent {\em Proof: }The main point is the easy observation that if $A$ has
$2$-components $A_1, \ldots, A_k$ then
$|N(A)|=|N(A_1)|+\ldots+|N(A_k)|$. Armed with this, we have
$\sum_{\mbox{$A \subseteq \cE$, $A$ small}} 2^{-|N(A)|}$
\begin{eqnarray*}
 & = & \sum_k
\sum_{\mbox{$A \subseteq \cE$ small with $2$-components
$A_1, \ldots, A_k$}} 2^{-|N(A_1)|-\ldots-|N(A_k)|} \\
& \leq & \sum_k \frac{1}{k!}\left(\sum_{\mbox{$A \subseteq \cE$
small, $2$-linked, $|A|\geq 1$}}
2^{-|N(A)|}\right)^k \\
& \leq & \exp\left\{\sum_{\mbox{$A \subseteq \cE$ small, $2$-linked,
$|A|\geq 1$}}
2^{-|N(A)|}\right\} \\
& = & \exp\left\{\sum_{\mbox{$A \subseteq \cE$, $|A|=1$}} 2^{-d} +
\sum_{\mbox{$A \subseteq \cE$ small, $2$-linked, $|A|\geq 2$}}
2^{-|N(A)|}\right\} \\
& = & \exp\left\{\frac{1}{2} + \sum_{\mbox{$A \subseteq \cE$ small,
$2$-linked, $|A|\geq 2$}} 2^{-|N(A)|}\right\}.
\end{eqnarray*}
The result will follow if we show that
\begin{equation} \label{sapsum}
\sum_{\mbox{$A \subseteq \cE$ small, $2$-linked, $|A|\geq 2$}}
2^{-|N(A)|} = o(1).
\end{equation}

Say that an $A \subseteq {\cal E}$ with $A$ small and $2$-linked is
{\em of type I} if $|N(A)|\leq d^4$ and of {\em of type II}
otherwise. We consider the portions of the sum in (\ref{sapsum})
corresponding to types I and II $A$'s separately. In each case we
partition the set of $A$'s according to the sizes of $[A]$ and
$N(A)$ and the first vertex $v \in N(A)$ (in some fixed ordering of
$\cO$).

If $A$ is of type I then by Lemma \ref{boundsondeltaAsmall}, $|[A]|
\leq O(|N(A)|/d)$ and so we may bound
\begin{equation} \label{Gagv_small}
|\cG(a,g,v)| \leq d2^{O(a\log d)}2^a = 2^{O(a\log d)} = 2^{O(g\log
d/d)}.
\end{equation}
The factor of $d$ comes from choosing a neighbour of $v$
to be a fixed vertex in $[A]$, the factor of $2^{(a\log d)}$ is from
Lemma \ref{Tree} and corresponds to a choice of $[A]$ (noting that
$[A]$ is $2$-linked) and the factor of $2^a$ comes from choosing $A$
as a subset of $[A]$. We therefore have
\begin{eqnarray}
\sum_{\mbox{$A$ of type I}} 2^{-|N(A)|} & \leq & \sum_{2 \leq a < g, ~g\leq d^4, ~v} |\cG(a,g,v)|2^{-g} \nonumber \\
& \leq & 2^d d^4 \sum_{g\geq d} 2^{O(g\log
d/d)-g} \label{explanation} \\
& = & e^{-\Omega(d)}. \label{typeIA}
\end{eqnarray}
The factor of $2^d$ in \eqref{explanation} comes from choosing $v$
while the factor of $d^4$ is for the choice of $a$.

For $A$ of type II (summing only over those values of $a$, $g$ and
$v$ with ${\cal G}(a,g,v) \neq \emptyset$ and with the inequalities
justified below)
\begin{eqnarray}
\sum_{\mbox{$A$ of type II}} 2^{-|N(A)|} & = & \sum_{a \leq g,~g \geq d^4,~v} |{\cal G}(a,g,v)| 2^{-g} \nonumber \\
  & \leq & \sum_{a \leq g,~g \geq d^4, ~v} 2^{-\Omega\left(\frac{t}{\log d}\right)} \label{usingmainapprox} \\
  &  \leq  & 2^{3d} 2^{-\Omega\left(\frac{d^{7/2}}{\log d}\right)} \label{msquared} \\
  &   =  & e^{-\Omega(d)}. \label{typeIIIAgpart}
\end{eqnarray}
Here (\ref{usingmainapprox}) is from Lemma \ref{lem-mainsaplemma}.
To see that the application is valid, note that
\eqref{lowerboundont} is given by Claim \ref{clm-hyper.iso}, and
that $\gD_2(Q_d)=2$. In (\ref{msquared}) we use the fact that there
are fewer than $2^d$ choices for each of $a$,  $g$ and $v$ and we
bound
$$
t = g\left(\frac{g-a}{g}\right) \geq
\Omega\left(\frac{g}{\sqrt{d}}\right) \geq \Omega(d^{7/2}),
$$
the first inequality following from Claim \ref{clm-hyper.iso} and
the second from the fact that $g \geq d^4$.

Combining (\ref{typeIA}) and (\ref{typeIIIAgpart}), we have
(\ref{sapsum}). \qed

\medskip

We may now swiftly put an upper bound on $|\cI(Q_d)|$. Indeed, it is
easy to check that for any $I \in \cI(Q_d)$ we have that $[I \cap
\cE]$ and $[I \cap \cO]$ have no edges between them, and so at least
one of $I \cap \cE$, $I \cap \cE$ is small. Noting that once we have
chosen a small $A \subseteq \cE$ there are $2^{2^{d-1}-|N(A)|}$ ways
to complete this choice to an independent set by selecting an
arbitrary subset of the non-neighbours of $A$ on $\cO$, by symmetry
it follows that
\begin{eqnarray*}
|\cI(Q_d)| & \leq & 2 \sum_{A \subseteq \cE} 2^{2^{d-1}-|N(A)|} \\
& = & 2(1+o(1))\sqrt{e}2^{2^{d-1}},
\end{eqnarray*}
the second inequality coming from Corollary \ref{cor-mainsaplemma}.

To obtain the matching lower bound and complete the proof of Theorem
\ref{thm-sapkor}, let $f(k)$ denote the number of subsets $S$ of
$\cE$ of size $k$ which satisfy the condition that $|N(S)|=kd$ (the
maximum possible; achievable only if the elements of $S$ have
pairwise disjoint neighbourhoods). Noting that for each $v \in \cE$
there are ${d \choose 2}$ vertices at distance $2$ from $v$ (and
exactly one at distance $0$), for $k \leq d$ (say) we have
\begin{eqnarray*}
f(k) & \geq & \frac{1}{k!}\prod_{j=0}^{k-1}
\left(2^{d-1}-(j-1)\left({d
\choose 2} + 1\right)\right) \\
& \geq & (1-o(d^{-1}))\frac{\left(2^{d-1}\right)^k}{k!}~~~~~\mbox{as
$d \rightarrow \infty$}.
\end{eqnarray*}
We will get our lower bound by considering those independent sets
which have $k$ vertices with non-overlapping neighbourhoods on one
side, and are arbitrary on the other side, for $k \leq d$. In the
first inequality in this count, the final term of $2^{2d^2}$ upper
bounds the overcount; the contribution from those sets which have at
most $d$ vertices from $\cE$ and at most $d$ from $\cO$.
\begin{eqnarray*}
|\cI(Q_d)| & \geq & 2\left(\sum_{k=0}^d f(k)2^{2^{d-1}-kd}\right) -
2^{2d^2}
\\
& \geq & 2^{2^{d-1}+1} \left(\sum_{k=0}^d
(1-o(d^{-1}))\frac{\left(2^{d-1}\right)^k}{k!} (1-o(d^{-1}))2^{-kd}\right)  - 2^{2d^2}\\
& \geq & 2(1-o(1))2^{2^{d-1}} \sum_{k=0}^d
\frac{\left(\frac{1}{2}\right)^k}{k!}
 \\
& \geq & 2(1-o(1))\sqrt{e}2^{2^{d-1}}.
\end{eqnarray*}
We have shown, as intended, that
$$
|\cI(Q_d)| \sim 2\sqrt{e}2^{2^{d-1}} ~~~~~~~\mbox{as $d \rightarrow
\infty$}.
$$

\section{Proof of Lemma \ref{lem-mainsaplemma}} \label{sec-saplemmaproof}

We bound $|{\cal G}(a,g,v)|$
using two notions of ``approximation''. These are introduced in
Section \ref{subsec-approximations}, and in this section we also
state the three ``approximation'' lemmas that we will use to obtain
the results discussed above. Section \ref{subsec-mainsaplemmaproofs}
gives the proof of Lemma \ref{lem-mainsaplemma},
modulo the approximation lemmata,
while Sections \ref{subsec-phiproof}, \ref{subsec-psiproof} and
\ref{subsec-reconstructionproof} are then devoted to the proofs of
the approximation lemmata.

\subsection{Approximation} \label{subsec-approximations}

The first notion of approximation depends on a parameter $\varphi$, $1
\leq \varphi \leq d-1$. Set
$$
G^{\varphi} = \{y \in G:d_{[A]}(y)>\varphi\}.
$$

\begin{defn}
A {\em $\varphi$-approximation} for $A \subseteq X$ is an $F' \subseteq
Y$ satisfying \beq{firstpprox1} G^{\varphi} \subseteq F' \subseteq G
\enq and \beq{firstpprox2} N(F') \supseteq [A] \enq
\end{defn}

The second depends on a parameter $\psi$, $1 \leq \psi \leq d-1$.
\begin{defn}
A {\em $\psi$-approximation} for $A \subseteq X$ is a pair $(F,S)
\in 2^Y \times 2^X$ satisfying \beq{secondapprox1} F \subseteq G,~~S
\supseteq [A], \enq \beq{secondapprox2} d_F(u) \geq d-\psi~~~\forall
u \in S \enq and \beq{secondapprox3} d_{X \setminus S}(v) \geq
d-\psi~~~\forall v \in Y \setminus F. \enq
\end{defn}

Before continuing, we note a property of $\psi$-approximations that
will be of use later.
\begin{lemma} \label{sleqf}
If $(F,S)$ is a $\psi$-approximation for $A \in {\cal G}$ then
\beq{boundingsbyf} |S| \leq |F| + 2t\psi/(d-\psi). \enq
\end{lemma}

\noindent {\em Proof: }Observe that $|\nabla(S,G)|$ is bounded above
by $d|F| + \psi|G \setminus F|$ and below by $d|[A]| + (d-\psi)|S
\setminus [A]| = d|S| - \psi|S \setminus [A]|$, giving
$$|S| \leq |F| + \psi|(G \setminus F) \cup (S \setminus [A])|/d,$$
and that each $u \in (G \setminus F) \cup (S \setminus [A])$
contributes at least $d-\psi$ edges to $\nabla(G,X\setminus [A])$, a
set of size $td$, giving
$$|(G \setminus F) \cup (S \setminus A)| \leq 2td/(d-\psi).$$
These two observations together give (\ref{boundingsbyf}). \qed

\medskip

In what follows we write ${\cal G}$ for ${\cal G}(a,g,v)$. We will
bound $|{\cal G}|$ by combining the following three lemmata.
\begin{lemma} \label{firstapprox}
For $g > d^4$ there is a family ${\cal V} = {\cal V}(\varphi) \subseteq
2^Y$ with \beq{howbigv} |{\cal V}| \leq \left\{
                  \begin{array}{ll}
                    2^{O(g \log^2 d/(\varphi d))+O(t \log^2 d/\varphi)} & \mbox{ if $t < O(g(d-\varphi)/(\varphi d))$ and} \\
                    2^{O(t \log^2 d/(d-\varphi))+O(t \log^2 d/\varphi)} & \mbox{ if $t > \Omega(g(d-\varphi)/(\varphi d))$}
                  \end{array}
                \right.
\enq such that each $A \in {\cal G}$ has a $\varphi$-approximation in
${\cal V}$.
\end{lemma}

\begin{lemma} \label{secondapprox}
For any $F' \in {\cal V}(\varphi)$ and $1 \leq \psi \leq d-1$ there is
a family ${\cal W}={\cal W}(F', \varphi, \psi) \subseteq 2^Y \times
2^X$ with \beq{howbigw} |{\cal W}| \leq 2^{O(td\log
d/((d-\varphi)\psi))+O(td\log d/((d-\psi)\psi))} \enq such that any $A
\in {\cal G}$ for which $F'$ is a $\varphi$-approximation has a
$\psi$-approximation in ${\cal W}$.
\end{lemma}

\begin{lemma} \label{countingas}
Given $1 \leq \psi \leq d-1$ and $\gamma > 0$, for each $(F,S) \in
2^Y \times 2^X$ that satisfies (\ref{boundingsbyf}) there are at
most \beq{howmanyas} \max \left\{
        2^{g - \gamma t}, 2^{g-t+O((t\psi/(d-\psi)+\gamma t)\log d)}
     \right\}
\enq $A$'s in ${\cal G}$ satisfying $F \subseteq G$ and $S \supseteq
[A]$.
\end{lemma}

\subsection{Derivation of Lemma \ref{lem-mainsaplemma}} \label{subsec-mainsaplemmaproofs}

For $t$ satisfying (\ref{lowerboundont}), we obtain
$$
|{\cal G}| < 2^{g-\Omega(t/\log d)}
$$
by taking $\gamma = c/\log d$, $\psi = c'd/\log d$ (for suitably
chosen constants $c, c'$) and, for example, $\varphi = d/2$.

\subsection{Proof of the $\varphi$-approximation lemma} \label{subsec-phiproof}

Our $\varphi$-approximation $F'=F'(A)$ for a particular $A \in {\cal
G}$ will consist of three pieces. The first of these is
$N(N_{[A]}(T_0))$ where $T_0$ is a small subset of $G$ for which
$N(N_{[A]}(T_0))$ contains most of $G^{\varphi}$ and for which
$\Omega:=\nabla(T_0,X\setminus [A])$ is also small. (A suitably
chosen random $T_0$ does both of these.) The second piece is
$T_0':=G^{\varphi} \setminus N(N_{[A]}(T_0))$. Setting
$L=N(N_{[A]}(T_0)) \cup T_0'$, we have $L \supseteq G^{\varphi}$. The
final piece, $T_1$, is a small subset of $G \setminus L$ whose
neigbourhood includes $[A] \setminus N(L)$ (we use Lemma
\ref{lovaszstein} to bound $|T_1|$). Clearly $F'=N(N_{[A]}(T_0))
\cup T_0' \cup T_1$ is a $\varphi$-approximation for $A$. We then take
${\cal U}$ to be the collection of $F'$'s that are produced in this
way as we run over all possible $A \in {\cal G}$.

To control $|{\cal U}|$ we observe that in this procedure each $F'$
is given by the quadruple $(T_0, T_0', T_1, \Omega)$, where $T:=T_0
\cup T_0' \cup T_1$ is a small $8$-linked subset of $G$ and $\Omega$
is a small subset $\nabla(T_0)$. Lemma \ref{Tree} bounds the number
of possible $T$'s, and direct calculations limit the number of
choices for  $T_0$, $T_1$ and $\Omega$ given $T_0$.

Fix $A \in {\cal G}$. Set $p=20\gD_2 \log d / (\varphi d)$.
\begin{claim} \label{clm-constructing.T_0}
There is a $T_0 \subseteq G$ such that $v \in T_0$ and \beq{t0}
|T_0|\leq 4gp \enq \beq{omega} |\nabla(T_0,X\setminus [A])| \leq
4tdp \enq and \beq{t0prime} |G^{\varphi} \setminus N(N_{[A]}(T_0))|
\leq 3g/d^{10}. \enq
\end{claim}

\noindent {\em Proof: }Construct a random subset $S$ of $G$ by
putting each $y \in G$ in $S$ with probability $p$, these choices
made independently. Clearly \beq{towardst0} {\bf E}(|S|)=gp \enq
\noindent and since $|\nabla(G,X\setminus [A])| = td$,
\beq{towardsomega} {\bf E}(|\nabla(S,X\setminus [A])|) = tdp. \enq

By
the co-degree condition, for $y \in G^{\varphi}$ we have
$$
|N(N_{[A]}(\{y\}))| \geq \frac{\varphi d}{2\gD_2}
$$
and so
\begin{eqnarray}
{\bf E}(|G^{\varphi} \setminus N(N_{[A]}(S))|) & = & \sum_{y \in
G^{\varphi}} {\bf P}(y \not \in
N(N_{[A]}(S))) \nonumber \\
& = & \sum_{y \in G^{\varphi}} {\bf P}(N(N_{[A]}(\{y\})) \cap T_0 =
\emptyset) \nonumber \\
& \leq & g(1-p)^{\frac{\varphi d}{2\gD_2}} \nonumber \\
& \leq & g/d^{10} \label{towardst0prime}
\end{eqnarray}

Combining (\ref{towardst0}), (\ref{towardsomega}) and
(\ref{towardst0prime}) and using Markov's inequality we find that
there is at least one $T_0^{initial} \subseteq G$ satisfying
$$
|T_0|\leq 3gp, ~~~ |\nabla(T_0,X\setminus [A])| \leq 3tdp
$$
and (\ref{t0prime}). Now note that $p> \Omega(\log d/d^2)$, so for
$g
> d^4$ we have (as usual, for sufficiently large $d$)
\begin{equation} \label{modified_T_0}
gp \geq  1 ~~~\mbox{and}~~~tdp  \geq d.
\end{equation}
Set $T_0 = T_0^{initial} \cup \{v\}$. By (\ref{modified_T_0}) $T_0$
satisfies (\ref{t0}) and (\ref{omega}), and it inherits
(\ref{t0prime}) from $T_0^{initial}$. \qed

Set $T_0' = G^{\varphi} \setminus N(N_{[A]}(T_0))$, $\Omega =
\nabla(T_0,X\setminus [A])$ and $L=N(N_{[A]}(T_0)) \cup T_0'$. Let
$T_1 \subseteq G \setminus L$ be a cover of minimum size of $[A]
\setminus N(L)$ in the graph induced by $(G \setminus L) \cup ([A]
\setminus N(L))$. Then $F'=L \cup T_1$ is a $\varphi$-approximation for
$A$.

Before estimating how many sets $F'$ might be produced in this way
as we run over all $A\in {\cal G}$, we make some observations about
the sets described above.

First, note that by Lemma \ref{nearbysets} $F'$ is $4$-linked ($A$
is $2$-linked, every $x\in A$ is at distance $1$ from $F'$ and every
$y \in F'$ is at distance $1$ from $A$) and so, again by Lemma
\ref{nearbysets}, $T=T_0 \cup T_0' \cup T_1$ is $8$-linked (every $y
\in T$ is at distance $2$ from something in $F'$ and every $y \in
F'$ is at distance $2$ from something in $T$).

By (\ref{t0}) we have $|T_0| \leq  O(g \log d/(\varphi d))$, by
(\ref{t0prime}) $|T_0'| \leq O(g/ d^{10})$, and by (\ref{omega})
$|\Omega| \leq O(t \log d/\varphi)$.

To bound $|T_1|$, note that $|G \setminus L| \leq td/(d-\varphi)$ (each
vertex in $G \setminus L$ is in $G \setminus G^{\varphi}$ and so
contributes at least $(d-\varphi)$ edges to $\nabla(G, X\setminus
[A])$, a set of size $td$), $d_{[A] \setminus N(L)}(u) \leq d$ for
each $u \in G \setminus L$, and $d_{G \setminus L}(v) = d$ for each
$v \in [A] \setminus N(L)$. So by Lemma \ref{lovaszstein}, $|T_1|
\leq (t/(d-\varphi)) (1+\ln d) = O(t\log d/(d-\varphi))$.

Combining these observations, we get that $T$ is an $8$-linked
subset of $Y$ with $|T|=O(g \log d/(\varphi d))$ (if $t <
O(g(d-\varphi)/(\varphi d))$) and $|T|=O(t\log d/(d-\varphi))$ (if $t >
\Omega(g(d-\varphi)/(\varphi d))$). We deal with these two cases
separately.

If $t < O(g(d-\varphi)/(\varphi d))$ we apply Lemma \ref{Tree} to find
that there are $2^{O(g \log^2 d/(\varphi d))}$ possible choices for $T$
(note that $v_0 \in T$). Once $T$ has been chosen, there are a
further $2^{O(g \log^2 d/(\varphi d))}$ choices for $T_0 \subseteq T$,
the same number of choices for $T_1 \subseteq T$ and $\sum_{i \leq
O(t \log d/\varphi)} {|\nabla(T_0)| \choose i} = 2^{O(t \log^2
d/\varphi)}$ choices for $\Omega$. So the total number of choices for
the quadruple $(T_0, T_0', T_1, \Omega)$ is
$$
2^{O(g \log^2 d/(\varphi d))+O(t \log^2 d/\varphi)}.
$$

If $t > \Omega(g(d-\varphi)/(\varphi d))$ we apply Lemma \ref{Tree} to
find that there are $2^{O(t \log^2 d/(d-\varphi))}$ possible choices
for $T$. Once $T$ has been chosen, there are a further $2^{O(t
\log^2 d/(d-\varphi))}$ choices for $T_0 \subseteq T$, the same number
of choices for $T_1 \subseteq T$ and (as before) $2^{O(t \log^2
d/\varphi)}$ choices for $\Omega$. So the total number of choices for
the quadruple $(T_0, T_0', T_1, \Omega)$ is
$$
2^{O(t \log^2 d/(d-\varphi))+O(t \log^2 d/\varphi)}.
$$

Once $T_0, T_1$ and $\Omega$ have been chosen, $F'$ is completely
determined, and (\ref{howbigv}) follows.

\subsection{Proof of the $\psi$-approximation lemma} \label{subsec-psiproof}

Fix a linear ordering $\ll$ of $V$. Given $A \in {\cal G}$ for which
$F'$ is a $\varphi$-approximation, we produce a $\psi$-approximation
$(F,S)$ for $A$ via the following algorithm.

\medskip

\noindent {\bf Step $1$: }If $\{u \in [A]: d_{G \setminus F'}(u) >
\psi \}
 \neq \emptyset$, pick the smallest (with respect to $\ll$) $u$ in this
set and update $F'$ by $F' \longleftarrow F' \cup N(u)$. Repeat this
until $\{u \in [A]: d_{G \setminus F'}(u) > \psi \} = \emptyset$.
Then set $F''=F'$ and $S''=\{u \in X:d_{F''}(u) \geq d-\psi\}$ and
go to Step $2$.

\medskip

\noindent {\bf Step $2$: }If $\{w \in Y \setminus G: d_{S''}(w) >
\psi \} \neq \emptyset$, pick the smallest (with respect to $\ll$)
$w$ in this set and update $S''$ by $S'' \longleftarrow S''
\setminus N(w)$. Repeat this until $\{w \in Y \setminus G:
d_{S''}(w)
> \psi \} = \emptyset$. Then set $S=S''$ and
$F=F''\cup \{w \in Y :d_S(w) > \psi \}$ and stop.

\medskip

\begin{claim} \label{algoanalysisoutput}
The output of this algorithm is a $\psi$-approximation for $A$.
\end{claim}

\noindent {\em Proof: }To see that $F \subseteq G$ and $S \supseteq
[A]$, first observe that $F'' \subseteq G$ (an immediate consequence
of $F' \subseteq G$ and the procedure in Step $1$) and that $S''
\supseteq [A]$ (or Step $1$ would not have terminated). We then have
$S \supseteq [A]$ since Step $2$ deletes from $S''$ only neighbours
of $Y \setminus G$, and $F \subseteq G$ since the vertices added to
$F''$ at the end of Step $2$ are all in $G$ (or Step $2$ would not
have terminated). This gives (\ref{secondapprox1}).

To verify (\ref{secondapprox2}) and (\ref{secondapprox3}), note that
$d_{F''}(u) \geq d-\psi~\forall u \in S''$ by definition,
$S\subseteq S''$, and $F \supseteq F''$, so that $d_{F}(u) \geq
d-\psi~\forall u \in S$, and if $w \in Y \setminus F$ then $d_S(w)
\leq \psi$ (by Step $2$), so that $d_{X \setminus S}(w) \geq
d-\psi~\forall w \in Y \setminus F$. \qed

\begin{claim} \label{algoanalysisnum}
The algorithm described above has at most
$$2^{O(td\log d/((d-\varphi)\psi))+O(td\log d/((d-\psi)\psi))}$$
outputs as the input runs over those $A \in {\cal G}$ for which $F'$
is a $\varphi$-approximation.
\end{claim}

\noindent Taking ${\cal W}$ to be the set of all possible outputs of
the algorithm, the lemma follows.

\medskip

\noindent {\em Proof of Claim \ref{algoanalysisnum}: } The output of
the algorithm is determined by the set of $u$'s whose neighbourhoods
are added to $F'$ in Step $1$, and the set of $w$'s whose
neighbourhoods are removed from $S''$ in Step $2$.

Initially, $|G \setminus F'| \leq td/(d-\varphi)$ (each vertex in $G
\setminus F'$ is in $G \setminus G^{\varphi}$ and so contributes at
least $(d-\varphi)$ edges to $\nabla(G, X\setminus [A])$, a set of size
$td$). Each iteration in Step $1$ removes at least $\psi$ vertices
from $G \setminus F'$ and so there can be at most
$td/((d-\varphi)\psi)$ iterations. The $u$'s in Step $1$ are all drawn
from $[A]$ and hence $N(F')$, a set of size at most $dg$. So the
total number of outputs for Step $1$ is at most
$$\sum_{i \leq td/((d-\varphi)\psi)} {dg \choose i} = 2^{O(td\log d/((d-\varphi)\psi))}.$$

We perform a similar analysis on Step $2$. Each $u \in S''\setminus
[A]$ contributes more than $d-\psi$ edges to $\nabla(G, X\setminus
[A])$, so initially $|S''\setminus [A]|\leq td/(d-\psi)$. Each $w$
used in Step $2$ reduces this by at least $\psi$, so there are at
most $td/((d-\psi)\psi)$ iterations. Each $w$ is drawn from
$N(S'')$, a set which is contained in the fourth neighbourhood of
$F'$ and so has size at most $d^4g$. So the total number of outputs
for Step $2$ is
$$2^{O(td\log d/((d-\psi)\psi))}.$$
The claim follows. \qed

\subsection{Proof of the reconstruction lemma}
\label{subsec-reconstructionproof}

Say that $S$ is {\em small} if $|S| < g - \gamma t$ and {\em large}
otherwise. We can obtain all $A \in {\cal G}$ for which $F \subseteq
G$ and $S \supseteq [A]$ as follows.

If $S$ is small, we specify of $A$ by picking a subset of $S$. If
$S$ is large, we first specify of $G$. Note that by
(\ref{boundingsbyf}) and the definition of large we have $|G
\setminus F| < 2t\psi/(d-\psi)+\gamma t$ and that $G \setminus F
\subseteq N(S) \setminus F$, so we specify $G$ by picking a subset
of $N(S) \setminus F$ of size at most $2t\psi/(d-\psi)+\gamma t$
(this is our choice of $G \setminus F)$. Then, noting that $[A]$ is
determined by $G$, we specify of $A$ by picking a subset of $[A]$.

This procedure produces all possible $A$'s (and lots more besides).
We now bound the total number of outputs.

If $S$ is small then the total number of possibilities for $A$ is at
most \beq{ssmall} 2^{g - \gamma t}. \enq

We have
$$|N(S)\setminus F| \leq d|S| \leq 3d^2g$$
so that if $S$ is large, the total number of possibilities for $|G
\setminus F|$ is at most
$$\sum_{i< 2t\psi/(d-\psi)+\gamma t}{3d^2g \choose i} \leq 2^{O((t\psi/(d-\psi)+\gamma t)\log d)}.$$
and so the total number of possibilities for $A$ is at most
\beq{slarge} 2^{g-t+O((t\psi/(d-\psi)+\gamma t)\log d)}. \enq

The lemma follows from (\ref{ssmall}) and (\ref{slarge}).

\end{document}